\documentclass[12pt,reqno]{amsart}
\usepackage{amssymb, amsfonts, amsbsy, latexsym}
\usepackage[english]{babel}

\newtheorem{Thm}{Theorem}[section]

\newtheorem{theorem}[Thm]{Theorem}

\newtheorem{lemma}[Thm]{Lemma}
\newtheorem{prop}[Thm]{Proposition}

\theoremstyle{remark}
\newtheorem{remark}[Thm]{Remark}

\DeclareMathOperator{\GL}{GL}
\DeclareMathOperator{\SL}{SL}
\DeclareMathOperator{\Gr}{Gr}

\DeclareMathOperator{\Sym}{Sym}

\def \A{\mathbb A}
\def \P{\mathbb P}
\def \G{\mathbb G}
\def \O{\mathcal O}
\def \V{\mathcal V}
\def \F{\mathcal F}
\def \W{\mathcal W}

\def \det{\text{det }}
\def \beq{\begin{eqnarray*}}
\def \eeq{\end{eqnarray*}}

\begin{document}

\title{Grassmannians and representations}
\author{Dan Edidin and Christopher A. Francisco}

\address{
Department of Mathematics,
University of Missouri,
Columbia, MO 65211}
\email{edidin@math.missouri.edu} 
\email{chrisf@math.missouri.edu}
\thanks{The first author was partially supported by N.S.A. grant
MDA904-03-1-0040}
\begin{abstract}
In this note we use Bott-Borel-Weil theory to
compute cohomology of interesting vector
bundles on sequences of Grassmannians.
\end{abstract}

\subjclass{Primary 14M15; Secondary 20G05}


\maketitle

\section{Introduction}\label{sec.intro}
For any positive integer $n$ the higher cohomology of
the line bundle $\O(n)$ on $\P^{m-1}$ vanishes. The dimension
of the space of global sections of this bundle is easily
calculated to be ${n+m \choose n}$ via the identification
of $H^0(\P^{m-1},\O(n))$ with the vector space of homogeneous
forms of degree $n$ in $m$ variables. 

If we view $\P^{m-1}$ as the space of lines in an $m$-dimensional
vector space $V$, then the line bundle $\O(n)$ is the $n$-th tensor
power of the dual of the tautological line subbundle $\O(-1)$. 
Generalizing to the
Grassmannian of $k$-planes we are led to a number of questions about
the cohomology of vector bundles on Grassmannians.

The most obvious question
is the following: Let $V$ be a vector space of dimension $m$. Compute the dimension of the space of sections $H^0(\Gr(k,V),
\V^{\otimes n})$ as a function of $k,m,n$, where $\V$ is the
dual of the tautological rank $k$ subbundle on the Grassmannian of $k$-planes in
$V$. Likewise, if we are interested in the dimension of linear systems
we could ask for the dimension of the linear system $H^0(\Gr(k,V),
(\det \V)^{\otimes n})$.

In this note we use Bott-Borel-Weil theory to answer these
questions. In particular we will compute for any (irreducible)
representation $W$ of $\GL_k$ the dimension $H^0(\Gr(k,m),{\mathcal W})$
where ${\mathcal W}$ is the vector bundle on $\Gr(k,m)$ whose
fiber at a point corresponding to a $k$-dimensional linear subspace $L$ 
is the dual vector space $W^*$ viewed as a $\GL(L)$-module. In addition
we explain when the higher cohomology vanishes.

Our motivation came from trying to understand
cohomology of vector bundles on classifying spaces.
It is also closely related to the study of representations
of $ind$-groups (\cite{twoivansandthewolf}). 
The classifying space, $B\G_m$,
for the multiplicative group $\G_m$,  is the infinite projective
space $\P^\infty$. 
The sequence of bundles $\{\O_{\P^m}(n)\}_{m=1}^\infty$ defines
a line bundle on $B\G_m$ whose
global sections is the limit of the countable sequence
of finite dimensional vector spaces $H^0(\P^m,\O(n))$.

In a similar way the classifying space $B\GL_k$ may be identified
with the infinite Grassmannian $\Gr(k,\infty)$ which is the limit
of the finite-dimensional Grassmannians $\Gr(k,m)$.
For a given $\GL_k$-module $W$, the associated
sequence of vector bundles on $\Gr(k,m)$ defines a vector
bundle on the classifying space $B\GL_k$. As we show below, the global
sections of these bundles on $Gr(k,m)$ have natural
$\GL_m$-module structures and the sequence of global sections  
$\{H^0(Gr(k,m),{\mathcal W})\}_{m=k}^\infty$ defines
a representation of
the $ind$-group $\GL_\infty$. When $W$ is irreducible this
is an irreducible representation of $\GL_\infty$ in the sense
of \cite{twoivansandthewolf}.

\section{Definitions and statement of results}
We work over an arbitrary algebraically closed field $K$.

For any $k$ and $m \geq k$ we may construct
the Grassmannian $\Gr(k,m)$ as follows: Let
$U_{k,m}\subset \A^{mk}$ be the open set parametrizing $m \times k$
matrices of rank $k$. There is a free action of $\GL_k$
on $U_{k,m}$ given by right matrix multiplication and
$\Gr(k,m) = U_{k,m}/\GL_k$.
The principal bundle $U_{k,m} \to \Gr(k,m)$ is the {\it frame
bundle}; i.e., the fiber over a point of $\Gr(k,m)$ corresponding
to a linear subspace $L$ is the set of all possible bases for $L$
in the vector space $K^m$. When $k =1$, this gives the familiar
construction of $\P^{m-1}$ as the quotient of $\A^m -\{0\}$ by
$\G_m$.

If $W$ is a (left) $\GL_k$-module then we obtain a vector bundle
on $\Gr(k,m)$ by taking the quotient $\W = U_{k,m}  \times_{\GL_k} W$
where $\GL_k$ acts diagonally on the product $U_{k,m} \times W$.
When $W$ is the defining
representation of $\GL_k$, then ${\mathcal W}$ is the dual of the 
tautological
rank $k$ subbundle of $\Gr(k,m)$. Likewise if $W$ is the determinant
character, then $\W$ is the line bundle which gives the Pl\"ucker
embedding.

Since $\GL_m$ acts on $U_{k,m}$ by left matrix multiplication
the quotient bundle $\W = U_{k,m} \times_{\GL_k} W$
is $\GL_m$-equivariant, so the cohomology of $\W$ is a (left) $\GL_m$-module.
The main result of this paper is the determination
of these cohomology modules. 

To state our theorem we recall some notation about $\GL_k$-modules.
Let $T_k \subset \GL_k$ be the group of diagonal matrices and let
$B_k$ the group of upper triangular matrices.  Any nonincreasing
sequence of integers $\lambda = (\lambda_1, \lambda_2,\ldots
\lambda_k)$ determines a $\GL_k$-module $V^{(k)}_\lambda$ which
is defined as follows:
The sequence of integers $\lambda$
determines a character of $T_k$ which we also call $\lambda$.
This character extends to a
representation of $B_k$. Let $L_\lambda$ be the
line bundle $\GL_k \times_{B_k} \lambda$, where
$B_k$ acts 
diagonally,\footnote{In \cite[p.393]{Fulton-Harris}
the line bundle whose global sections is the highest
weight module $V_\lambda$ is denoted $L_{-\lambda}$.}
and set $V_\lambda^{(k)} = 
H^0(\GL_k/B_k,L_\lambda)$. 
When the characteristic of $K$ is 0 then
$V^{(k)}_\lambda$ is irreducible. In any characteristic,
the groups $H^i(GL_k/B_k, L_{\lambda})$ vanish for
$i > 0$ \cite{Demazure}.

\begin{remark} With this notation $V_\lambda$ and $V_{\lambda'}$ restrict
to the same representation of $\SL_k$ if and only if
$\lambda_l - \lambda'_l$ is constant (i.e. independent of $l$).
In this case $V_{\lambda} = V_{\lambda'} \otimes D^{\otimes r}$
where $r = \lambda_l - \lambda'_l$ and $D$ is the determinant
character of $\GL_k$.
\end{remark}

Let $\lambda = (\lambda_1, \ldots , \lambda_k)$ be a sequence
of non-increasing integers and let $V_\lambda^{(k)}$ be the associated
representation of $\GL_k$. 
Let $\V_\lambda$ be the corresponding
vector bundle on $\Gr(k,m)$.

\begin{theorem} \label{thm.main}
(a) If $\lambda_k \geq 0$, then  
$H^0(\Gr(k,m), \V_\lambda) = V_\lambda^{(m)}$,
and $H^i(\Gr(k,m), \V_\lambda) = 0$ for $i>0$.
Here $V_\lambda^{(m)}$ denotes the $\GL_m$-module
with highest weight $(\lambda_1, \ldots , \ldots \lambda_k, 0, \ldots , 0)$.\\

(b) If $\lambda_k < 0$, then for sufficiently large $m$, 
$H^i(\Gr(k,m), \V_\lambda) = 0$ for all $i$.
\end{theorem}

\section{Bott-Borel-Weil theory for parabolics}
The material here is well known, but we do not know a reference with
algebraic proofs.

Let $H$ be an affine algebraic group and let $\pi\colon X \to Y$ be an
$H$-principal bundle.\footnote{This means \cite[Definition
0.10]{GIT} that $\pi$ is flat, $X \to Y$ is a geometric quotient and
that $X \times_Y X$ is $H$-equivariantly isomorphic to $H \times
X$. This definition implies that the action of $H$ on $X$ is free.}
If $\F$ is an $H$-equivariant $\O_X$-module, set
$\F^H$ to be subsheaf of invariant sections of $\pi_*\F$. 
Given an $\O_Y$-module ${\mathcal G}$, $\pi^*{\mathcal G} = 
\pi^{-1}{\mathcal G}
\otimes_{\pi^{-1}\O_Y} \O_X$ has a natural $H$-action given
by the action of $H$ on $\O_X$. Because 
$H$ acts freely, the functors ${\mathcal G} \mapsto \pi^*{\mathcal G}$
and $\F \mapsto \F^H$ are inverse to each other.

Let $X' \stackrel{\pi'} \to Y'$ and $X \to Y$ be principal
$H$-bundles. Let $q \colon X' \to X$ be a $H$-equivariant.
Then there is an induced map $p\colon Y' \to Y$ such that the
diagram
\begin{equation} \label{diag.cartesian}
\begin{array}{ccc}
X' & \stackrel{q} \to & X\\
\pi' \downarrow  & & \pi \downarrow\\
Y' & \stackrel{p} \to & Y
\end{array}
\end{equation}
commutes.

\begin{lemma} \label{lem.commute}
Diagram \eqref{diag.cartesian} is cartesian.
\end{lemma}
\begin{proof}
Since the diagram commutes there is an $H$-equivariant map
$X' \to Y' \times_Y X$. Now $X' \to Y'$ and $Y' \times_Y X \to Y'$
are both principal bundles over $Y$.  After base change by
a flat surjective morphism $X' \to Y'$ both bundles
are trivialized. Thus the map of $Y'$-schemes, 
$X' \to Y' \times_Y X$ is an isomorphism after flat surjective
base change. Therefore it is an isomorphism.
\end{proof}

\begin{lemma} \label{lem.invcommute}
Let
$\begin{array}{ccc}
X' & \stackrel{q} \to & X\\
\pi' \downarrow  & & \pi \downarrow\\
Y' & \stackrel{p} \to & Y
\end{array}$
be a commutative (hence cartesian) diagram of principal
$H$-bundles as above. Then for any $H$-equivariant $\O_{X'}$-module $\F$
there are natural isomorphism of $\O_Y$-modules
$(R^iq_* \F)^H = R^ip_*(\F^H)$ for any $i\geq 0$.
\end{lemma}
\begin{proof}

Because $\F = \pi^{\prime *}\F^H$, and
cohomology commutes with base change \cite[Proposition 9.3]{Hartshorne}, 
there are natural isomorphisms
$R^i q_* \F = \pi^* R^ip_*\F^H$. Pushing forward by
$\pi$ and taking $H$-invariants yields the desired isomorphism.
\end{proof}

Let $G$ be an algebraic group, $B \subset G$ a Borel
subgroup and $P \supset B$ a parabolic subgroup.
Let $p \colon G/B \to G/P$ be the projection.
\begin{prop} \label{prop.parabolic}
If $\lambda$ is a $B$-module, then $G \times_P H^i(P/B, P \times_B \lambda)
= R^ip_*(G \times_B \lambda)$ as $G$-equivariant bundles on $G/P$.
\end{prop}

\begin{proof}
The map $G \times P/B \to G/B$ given by $(g,pB) \mapsto gpB$
is a $P$-principal bundle, where $P$ acts freely on $G \times P/B$
by the formula $q\cdot(g,pB) = (gq, q^{-1}pB)$ where $q \in P$.

Thus we have a commutative diagram with the vertical
arrows proper and the horizontal arrows quotient maps by the free action
of $B$.
$$\begin{array}{ccc}
G \times P/B & \stackrel{q} \to & G\\
\downarrow & & \downarrow\\
G/B &\stackrel{p} \to & G/P
\end{array}$$
Consider the $P$-equivariant line bundle $L =G \times (P \times_B \lambda)$
on $G \times P/B$.
Since $q$ is a projection,
$R^iq_*L = G \times H^i(P/B, P \times_B \lambda)$.
On the other hand $L^P = G \times_B \lambda$ so
so by
Lemma \ref{lem.invcommute},
$G \times_P H^i(P/B, P \times_B \lambda)$ is naturally isomorphic
to $R^ip_*(G \times_B \lambda)$. Since this
isomorphism is natural, it is equivariant for the $G$-actions
on these bundles.
\end{proof}
\section{Proof of Theorem \ref{thm.main}}
Let $G = \GL_m$, and let $B_m$ be the Borel subgroup 
of upper triangular matrices. Let $P$ be the parabolic
subgroup of matrices of the form
$\left(\begin{array}{cc}
A_k & U \\
0 & A_{m-k}
\end{array}\right)$ with $A_k \in \GL_{k}$, $A_{m-k} \in \GL_{m -k}$ and $U$ an
arbitrary $k \times (m-k)$ matrix.  The subgroup $N$ of matrices of
the form $\left(\begin{array}{cc} I & U \\ 0 &
A_{m-k}\end{array}\right)$ is normal and $P/N = \GL_k$.

\begin{lemma} \label{lem.framebundle}
$U_{k,m}$ is isomorphic to the homogeneous space $G/N$.
\end{lemma}
\begin{proof}
Any $m \times k$ matrix of maximal rank can be viewed as the first $k$
columns of a nonsingular $m \times m$ matrix.  Since $\GL_m$ acts
transitively on itself by left multiplication, it also acts
transitively on the space $U_{k,m}$ of $m \times k$ matrices of
maximal rank. The stabilizer of the $m \times k$ matrix $\left(\begin{array}{c} I_k \\ 0\end{array} \right)$ is
the subgroup $N \subset \GL_k$.
\end{proof}

Using Lemma \ref{lem.framebundle} we obtain the well
known identification of $\Gr(k,m)$ $= U_{k,m}/\GL_k = G/P$. Likewise for
any $\GL_k$-module,
the vector bundle $U_{k,m} \times_{\GL_k} V$ is identified
with $G \times_P V$ where $V$ is made into a $P$-module via
the surjective map $P \to P/N = \GL_k$.

Since $P$ is a parabolic subgroup containing $B$
we have a proper map $G/B \stackrel{p} \to G/P$. 
Let $L_\lambda$ be the line bundle on $G/B$
corresponding to the weight $\lambda = (\lambda_1, \ldots , \lambda_k,
0, \ldots , 0)$. Then by Bott-Borel-Weil $H^0(G/B, L_\lambda) = 
V^{(m)}_\lambda$.

\begin{lemma} \label{lem.pushforward}
$p_*(L_\lambda)= \V_\lambda$ as $\GL_m$ equivariant bundles
on $\Gr(k,m)$.
\end{lemma}
\begin{proof}
By Proposition \ref{prop.parabolic}, $p_*(L_\lambda) = G \times_P H^0(P/B,
P \times_B L_\lambda)$. Thus to prove the Lemma
we must compute $H^0(P/B, P \times_B L_\lambda)$.

Since $P$ is a unipotent extension of $GL_k \times GL_{m-k}$,
the homogeneous space $P/B$ is isomorphic to the
product $GL_k/B_k \times GL_{m-k}/B_{m-k}$,
where $B_k = B \cap (\GL_k \times I_{m-k})$ and 
$B_{m-k} = B \cap (I_k \times \GL_{m-k})$. Thus 
\begin{eqnarray}
H^0(P/B, P \times_B \lambda) & = & H^0(\GL_k/B_k, \GL_k \times_{B_k} 
(\lambda|B_k)) \otimes \nonumber \\
& & H^0(\GL_{m-k}/B_{m-k}, \GL_{m-k} \times_{B_{m-k}} (\lambda|B_{m-k})). \nonumber
\end{eqnarray}

Since $\lambda = (\lambda_1, \ldots , \lambda_k, 0, \ldots ,0)$,
its restriction $B_k$ is the highest weight vector $(\lambda_1, \ldots
, \lambda_k)$, and its restriction to $B_{m-k}$ is trivial.
\end{proof}

Part (a) of the theorem now follows easily from Lemma \ref{lem.pushforward}.
First observe that $H^0(\Gr(k,m),\V_\lambda) = H^0(G/P,p_*L_\lambda)
= H^0(G/B, L_\lambda)$.
If all of the $\lambda_k$'s are nonnegative then 
$H^0(G/B, L_\lambda) = V^{(m)}_{\lambda}$ and 
$H^i(G/B, L_\lambda) = 0$ for $i > 0$.
It follows from the Leray spectral
sequence that
that $H^i(G/P, p_*L_\lambda) = 0$ as well.

We now prove (b). Suppose that $\lambda_k < 0$. It suffices
to show that, for $m$ sufficiently large, that $H^i(G/B, G \times_B
\lambda) = 0$ where is the character of the maximal torus
with weight $\lambda = (\lambda_1,\ldots , \lambda_k, 0, \ldots , 0)$.

Let $\{0\}=V_0 \subset V_1 \ldots \subset V_m$ be the
flag fixed by the Borel subgroup of upper triangular
matrices $B_m \subset \GL_m$. For any $l < m$
let $P_l$ be the parabolic subgroup fixing the subflag 
$V_0 \subset V_1 \subset \ldots \subset V_l$.  Then $P_{l} \subset 
P_{l-1}$ and the fibers of the map $p_l\colon G/P_{l} \to G/P_{l-1}$
are isomorphic to the projective space $\P^{m-l}$.
In this way the flag variety $G/B = G/P_{m-1}$ is realized
as being at the top of the tower of projective bundles
$$G/B \stackrel{p_{m-1}} \to G/P_{m-2} \to \ldots \to G/P_2 \stackrel{p_2}
\to G/P_1 = \P^{m-1}.$$

The group $P_l$ is isomorphic to the product $B_l \times \GL_{m-l}$
where $B_l$ is the Borel group of upper triangular matrices in
$\GL_l$. Thus, irreducible representations of $P_l$ are of
the form $\chi \times W$ where $\chi$ is a character of
the maximal torus in $B_l$ and $W$ is an irreducible representation
of $GL_{m-l}$.
If $\chi = (a_1, \ldots , a_l)$ is a character of the maximal
torus of $B_l$ 
then the same argument used in the proof Proposition \ref{prop.parabolic}
shows that $R^ip_{l*} (G \times_{P_l} \chi) = 
G \times_{P_{l-1}} H^i(P_{l-1}/P_{l}, P_{l-1} \times_{P_l} \chi)$.
Under the identification $P_{l-1}/P_{l} = \P^{m-l +1}$ the 
line bundle $P_{l-1} \times_{P_l} \chi$ corresponds to 
$\O_{\P^{m-l+1}}(a_l)$. 

If $\lambda$ is a character of $T_m$ with 
$\lambda_{k+1} = \ldots \lambda_m = 0$ then
we may view $\lambda$ as $P_l$-module for $l \geq k$.
The argument of the previous paragraph
implies that  $H^i(G/B, G \times_B \lambda) = H^i(G/P_k, G \times_{P_k}
\lambda)$. If $\lambda_k < 0$ then, for all $j$, 
$H^j(\P^{m-k}, \O_{\P^{m-k}}(\lambda_k)) = 0$ whenever $m \geq -\lambda_k +k$.
Thus $R^jp_{k*}(G/\times_{P_k} \lambda) = 0$ for all $j$ as long
as $m$ is sufficiently large.
Therefore, by the Leray spectral sequence, it follows that
$H^i(G/B, L_\lambda) = 0$ for all $i$.

\section{Some explicit dimension computations}\label{sec.formulas}
We conclude by using Theorem~\ref{thm.main} to compute the dimensions
of $H^0(\Gr(k,m),\V_\lambda)$ for some particular $\V_\lambda$
involving symmetric powers and the determinant. To do this, we combine
Theorem~\ref{thm.main} with \cite[Theorem 6.3]{Fulton-Harris} to
obtain a convenient formula for the dimensions.

Let $\lambda =
(\lambda_1, \ldots , \lambda_k)$ be a sequence of nonincreasing
integers, and let $V_\lambda^{(k)}$ be the associated representation
of $\GL_k$ with $\V_\lambda$ the corresponding vector bundle on
$\Gr(k,m)$.

\begin{lemma} \label{lem.dimensions}
If all $\lambda_i \ge 0$, then \[ \dim H^0(\Gr(k,m),\V_\lambda) = \prod_{i < j} \frac{\lambda_i - \lambda_j +j-i}{j-i},\] where $j$ runs from 2 to $m$, and $\lambda_{k+1}= \cdots = \lambda_m = 0$. 
\end{lemma}
 
The first proposition gives a symmetry result about the dimension of a
power of the determinant. Throughout the rest of this section, let $V$
denote the standard $k$-dimensional representation of $\GL_k$ and $\V$
the corresponding tautological rank $k$ vector bundle on $\Gr(k,m)$.  

\begin{prop} \label{prop.symmetry}
For all $0 < k \le m$ and $l > 0$, 
\[ \dim H^0(\Gr(k,m),(\det \V)^l) = \dim H^0(\Gr(k,k+l),(\det \V)^{m-k})\] 
\[ = \prod_{j=k+1}^m \frac{{l+j-1 \choose l}}{{l+j-k-1 \choose l}}. \]
\end{prop}

\begin{proof}
The partition $\lambda$ associated to $(\det \V)^l$ is $\lambda=(l,\dots,l,0,\dots,0)$. Hence by Lemma~\ref{lem.dimensions}, \[ \dim H^0(\Gr(k,m),(\det \V)^l) = \prod_{j=k+1}^m \prod_{i=1}^k \frac{l+j-i}{j-i} \] 
\[ = \prod_{j=k+1}^m \frac{l!(l+j-1)!(j-k-1)!}{l!(j-1)!(l+j-k-1)!} = \prod_{j=k+1}^m \frac{{l+j-1 \choose l}}{{l+j-k-1 \choose l}} =  \prod_{j=1}^{m-k} \frac{{l+j+k-1 \choose l}}{{l+j-1 \choose l}}, \]  
\[ \dim H^0(\Gr(k,k+l),(\det \V)^{m-k}) = \prod_{j=k+1}^{k+l} \prod_{i=1}^k \frac{m-k+j-i}{j-i} \] 
\[ = \prod_{j=k+1}^m \frac{(m-k)!(m-k+j-1)!(j-k-1)!}{(m-k)!(m+j-1)!(j-1)!} = \prod_{j=k+1}^{k+l} \frac{{m-k+j-1 \choose m-k}}{{m-k+j-k-1 \choose m-k}} \] 
\[ = \prod_{j=1}^l \frac{{m+j-1 \choose m-k}}{{m+j-k-1 \choose m-k}}.\]

To show that the two dimensions are equal, we induct on $l$. When $l=1$, it is easy to compute that both dimensions are ${m \choose k}$. Suppose now that $s > 1$, and we have proven that the dimensions are the same for $l=s-1$. The induction hypothesis gives us that 
\[ \prod_{j=1}^{s-1} \frac{{m+j-1 \choose m-k}}{{m+j-k-1 \choose m-k}} = \prod_{j=1}^{m-k} \frac{{s+j+k-2 \choose s-1}}{{s+j-2 \choose s-1}}.\] 
Hence, letting $l=s$, \[ \prod_{j=1}^{s} \frac{{m+j-1 \choose m-k}}{{m+j-k-1 \choose m-k}} = \frac{{m+s-1 \choose m-k}}{{m+s-k-1 \choose m-k}} \prod_{j=1}^{m-k} \frac{{s+j+k-2 \choose s-1}}{{s+j-2 \choose s-1}} \] 
\[ = \frac{{m+s-1 \choose m-k}}{{m+s-k-1 \choose m-k}} \prod_{j=1}^{m-k} \frac{s+j-1}{s+j+k-1} \prod_{j=1}^{m-k} \frac{{s+j+k-1 \choose s}}{{s+j-1 \choose s}}.\]

The rightmost product is what we want, so it suffices to show that the
product of the first two factors is one. Note that
\[ \prod_{j=1}^{m-k} \frac{s+j-1}{s+j+k-1} = \frac{\frac{(s+m-k-1)!}{(s-1)!}}{\frac{(s+m-1)!}{(s+k-1)!}} = \frac{{m+s-k+1 \choose m-k}}{{m+s-1 \choose m-k}},\] which completes the proof.
\end{proof}

\begin{remark} Assume that the characteristic of the ground
field is $0$.  If we fix $k$ and $l$ and let $m$ go to infinity, then
sequence of irreducible $\GL_m$ representations $\{V^{(m)}_{(l,\ldots
, l)}\}$ determines an irreducible representation of the $ind$-group
$\GL_\infty$. Proposition \ref{prop.symmetry} implies that the
dimensions of the finite dimensional representations $V^{(m)}_{(l, l,
\ldots ,l)}$ are the same as the dimensions of the irreducible
$\GL_{k+l}$-modules $V^{(l+k)}_{(m-k,\ldots , m-k)}$.
\end{remark}

\begin{remark}
As noted in \cite{twoivansandthewolf} the infinite Grassmannian
$\Gr(k,\infty)$ is an $ind$-projective variety.  For fixed $l$, the
sequence of $\GL_m$-modules $\{\Sym^l(\wedge^k V_m)/V^{(m)}_{(l,\ldots
, l)}\}$ is the $ind$-representation of $\GL_\infty$ corresponding to
the degree $l$ component of the ideal of Pl\"ucker relations for the
embedding of $\Gr(k,\infty)$ into the infinite projective space given
by the sequence of projective spaces $\{\P(\wedge^k V_m)\}$. (Here
$V_m$ denotes the standard representation of $\GL_m$.)

Proposition \ref{prop.symmetry} allows us to compute the 
sequence of dimensions of the Pl\"ucker relations in degree
$l$. In particular the dimension of the $\GL_m$-module of 
Pl\"ucker relations in degree $l$ is 
$${ {m \choose k} + l \choose l} - \prod_{j=k+1}^m \frac{{l+j-1
\choose l}}{{l+j-k-1 \choose l}}.$$  
\end{remark}

\begin{prop} \label{prop.sym}
For all $r \ge 0$ and $0 < k \le m$, \[ \dim H^0(\Gr(k,m),\Sym^r \V) =
{r+m-1 \choose r}. \]
\end{prop}

\begin{proof}
Use Lemma~\ref{lem.dimensions} and the partition $(r,0,\dots,0)$.
\end{proof}

\begin{prop} \label{prop.symdet}
For all $r \ge 0$, $l \ge 0$, and $0 < k \le m$, \[ \dim H^0(\Gr(k,m),\Sym^r \V \otimes (\det \V)^l) \]
\[= \frac{{r+l+m-1 \choose r}}{{r+l+k-1 \choose l}} {l+m-1 \choose l} \prod_{j=k+1}^m \frac{{l+j-2 \choose l}}{{l+j-k-1 \choose l}}. \]
\end{prop}

As we would expect, when $l=0$, we get $\dim H^0(\Gr(k,m),\Sym^r \V)$, and when $r=0$, we have (after a bit of manipulation) the dimension of $H^0(\Gr(k,m),(\det \V)^l)$. Also, the product over $j$ is exactly $\dim H^0(\Gr(k-1,m-1),(\det \V)^{l})$ (again easy to see after a little algebra). 

\begin{proof}
We use the partition $\lambda=(r+l,l,\dots,l,0,\dots,0)$, where $\lambda_{j} = 0$ for all $k+1 \le j \le m$. By Lemma~\ref{lem.dimensions}, the dimension is \[ \prod_{j=2}^k \frac{r+j-1}{j-1} \prod_{j=k+1}^m \frac{r+l+j-1}{j-1} \prod_{j=k+1}^m \prod_{i=2}^k \frac{l+j-i}{j-i}.\]

The first product is ${r+k-1 \choose r}$. Consider the second
product. We have
\[ \prod_{j=k+1}^m \frac{r+l+j-1}{j-1} = \frac{\frac{(r+l+m-1)!}{(r+l+k-1)!}}{\frac{(m-1)!}{(k-1)!}} = \frac{(r+l+m-1)!(k-1)!}{(r+l+k-1)!(m-1)!} \] 
\[ = \frac{(r+l+m-1)!}{r!(l+m-1)!} \frac{r!(l+m-1)!(k-1)!}{(r+l+k-1)!(m-1)!} \frac{l!(r+k-1)!}{l!(r+k-1)!} \] 
\[ = {r+l+m-1 \choose r} {r+l+k-1 \choose l}^{-1} {l+m-1 \choose l} {r+k-1 \choose r}^{-1}. \]

The last factor is \[ \prod_{j=k+1}^m \prod_{i=2}^k \frac{l+j-i}{j-i} = \prod_{j=k+1}^m \frac{\frac{(r+l+m-1)!}{(r+l+k-1)!}}{\frac{(m-1)!}{(k-1)!}} = \prod_{j=k+1}^m \frac{(l+j-2)!(j-k-1)!}{(l+j-k-1)!(j-2)!}.\] Multiplying the top and bottom by $l!$ yields \[ \prod_{j=k+1}^m \frac{{l+j-2 \choose l}}{{l+j-k-1 \choose l}}. \] 

Finally, multiplying all three products together completes the proof.
\end{proof}

\begin{remark}
Using Theorem~\ref{thm.main} it is easy to show that if $d \leq k$ then
$H^0(\Gr(k,m), \V^{\otimes d}) = V_m^{\otimes d}$ where
$V_m$ is the standard representation of $\GL_m$.
In particular $\dim H^0(\Gr(k,m), \V^{\otimes d}) = m^d$.
An interesting question is to determine whether 
there are relatively simple formulas
for the dimension of
$H^0(\Gr(k,m),\V^{\otimes d} \otimes (\det \V)^l)$.
If $l \geq  0$ then Theorem \ref{thm.main} gives the summation formula
$$\dim H^0(\Gr(k,m),\V^{\otimes d} \otimes (\det \V)^l) =
\sum_{\lambda} m_\lambda\; \dim V^{(m)}_{\lambda +l}$$
where the sum is over all $k$-element partitions $\lambda$ of $d$,
$m_\lambda$ is the product of the hook lengths
in the tableau corresponding to the partition 
$(\lambda_1, \ldots \lambda_k)$, and 
$V^{(m)}_{\lambda +l}$ refers to the $\GL_m$-module with
highest weight $(\lambda_1 + l, \ldots \lambda_k +l, 0, \ldots 0)$.
Unfortunately, we have do not know how to simplify this sum.
Other summation formulas can also be obtained using the Riemann-Roch
or localization theorems. 
When $-d < l < 0$,  there is the added difficulty
that some, but not all, of the sequences $(\lambda_1 + l, \lambda_2
+l, \ldots , \lambda_k +l, 0, \ldots 0)$ have negative entries,
so the corresponding line bundles on $\GL_m/B_m$ have no
cohomology for sufficiently large $m$.
\end{remark}

\end{document}